\documentclass[pdflatex,sn-mathphys-num]{sn-jnl}

\usepackage{amsmath,amssymb,amsfonts}
\usepackage{booktabs}
\usepackage{tabularx}
\usepackage{array}
\usepackage{enumitem}
\usepackage{tikz}
\usetikzlibrary{arrows.meta,positioning,calc}

\newcolumntype{L}[1]{>{\raggedright\arraybackslash}p{#1}}
\newcolumntype{Y}{>{\raggedright\arraybackslash}X}

\DeclareMathOperator{\spann}{span}
\DeclareMathOperator{\ad}{ad}
\DeclareMathOperator{\Ad}{Ad}
\DeclareMathOperator{\End}{End}
\DeclareMathOperator{\Hom}{Hom}
\DeclareMathOperator{\Lie}{Lie}
\DeclareMathOperator{\SO}{SO}
\DeclareMathOperator{\id}{id}

\theoremstyle{thmstyleone}
\newtheorem{theorem}{Theorem}[section]
\newtheorem{proposition}[theorem]{Proposition}
\newtheorem{lemma}[theorem]{Lemma}
\newtheorem{corollary}[theorem]{Corollary}

\theoremstyle{thmstylethree}
\newtheorem{definition}[theorem]{Definition}
\newtheorem{example}[theorem]{Example}

\theoremstyle{thmstyletwo}
\newtheorem{remark}[theorem]{Remark}

\begin{document}

\articletype{Research Article}

\title[Bialgebra structures on flat Lie algebras]{Bialgebra Structures on Flat Lie Algebras and Their Poisson--Lie Groups}

\author*[1]{\fnm{Amine} \sur{Bahayou}}
\email{amine.bahayou@gmail.com}

\affil[1]{\orgdiv{Department of Mathematics}, \orgname{Kasdi Merbah University}, \orgaddress{\city{Ouargla}, \postcode{30000}, \country{Algeria}}}

\abstract{%
We study Lie bialgebra structures on \emph{flat metric Lie algebras}, that is, Lie algebras $(\mathfrak{g},\langle\cdot,\cdot\rangle)$ whose associated left-invariant Riemannian metric on the simply connected Lie group $G$ has zero curvature. By Milnor's structure theorem, such $\mathfrak{g}$ splits orthogonally as
\[
\mathfrak{g}=\mathfrak{a}\oplus\mathfrak{u},\qquad \mathfrak{u}=[\mathfrak{g},\mathfrak{g}]\ \text{abelian and even dimensional},\quad\mathfrak{a}:=\mathfrak{s}\oplus\mathfrak{z},
\]
where $\mathfrak{z}$ is the center and $\mathfrak{s}$ is an abelian subalgebra that acts on $\mathfrak{u}$ by commuting infinitesimal rotations; this yields a decomposition of $\mathfrak{u}$ into $2$-dimensional weight planes $P_\ell$. Under a generic \emph{nondegeneracy} (nonresonance) condition on the weights, we establish a normal form for Lie-bialgebra $1$-cocycles $\xi\colon\mathfrak{g}\to \wedge^2\mathfrak{g}$: each $\xi$ admits a decomposition $\xi=\ad r+R$, where $\ad r$ is a coboundary and $R$ is a normalized cocycle with tightly controlled components. Using the Big Bracket (Maurer--Cartan) formalism together with the rotation geometry of the weight planes, we split the co-Jacobi condition into two independent equations: a reduced co-Jacobi equation $\{R,R\}=0$ for the normalized cocycle, and an invariant-trivector condition $[r,r]+2\{r,R\}\in(\wedge^3\mathfrak{g})^{\mathfrak{g}}$ for the coupling term. We then describe the quasi-triangular (classical Yang--Baxter) locus via invariant Schouten squares. Finally, we integrate $\xi$ to explicit multiplicative Poisson tensors on $G$, producing concrete families of flat Poisson--Lie groups with polynomial formulas along the abelian normal subgroup $\exp(\mathfrak{z}\oplus\mathfrak{u})$.
}


\maketitle

\section{Introduction}

Left-invariant flat Riemannian metrics on Lie groups sit at a natural crossroads of Lie theory, transformation groups, and invariant geometric structures. If a simply connected Lie group $G$ carries a left-invariant flat metric, then $G$ acts simply transitively by isometries on Euclidean space. In particular, the Lie algebra $\mathfrak{g}=\Lie(G)$ is highly constrained: it is two-step solvable and its commutator ideal is abelian. From the viewpoint of transformation groups, such groups provide canonical solvable models for Euclidean geometries. Moreover, many invariant tensors on $G$ (connections, Poisson structures, \emph{etc.}) can be studied effectively by combining the representation theory of the isotropy action with Lie-algebra cohomology.

In Poisson--Lie theory, a multiplicative Poisson tensor on $G$ is equivalent to a Lie bialgebra structure on $\mathfrak{g}$. Consequently, the description of Lie bialgebra structures on flat metric Lie algebras produces an explicit and geometrically meaningful class of Poisson--Lie groups adapted to Euclidean-type symmetry. The present paper develops a representation-theoretic approach tailored to the structure of flat metric Lie algebras, guided by two complementary goals:
\begin{enumerate}
\item \textbf{Structural control of the infinitesimal data.} A Lie bialgebra structure is a $1$-cocycle
\(
\xi:\mathfrak{g}\to\wedge^2\mathfrak{g}
\)
whose transpose endows $\mathfrak{g}^*$ with a Lie bracket. For flat $\mathfrak{g}$, the adjoint action admits a decomposition into weight planes; this can be exploited to place cocycles in a normal form and to separate potential obstructions weight by weight.
\item \textbf{Geometric realization on the group.} We integrate the resulting infinitesimal cocycles to explicit group $1$-cocycles
\(
\Pi:G\to\wedge^2\mathfrak{g},
\)
and hence to concrete multiplicative Poisson tensors on $G$, emphasizing formulas compatible with the semidirect-product structure of flat groups.
\end{enumerate}

The study of Lie bialgebra structures on flat Lie algebras is motivated, among other reasons, by the following considerations:
\begin{itemize}
\item \textbf{Links with deformation quantization.} Poisson--Lie groups play a central role in deformation quantization by bridging classical and noncommutative geometry. In particular, Hawkins \cite{Ha1,Ha2} related flatness-type conditions to deformation phenomena arising in formal quantization of Poisson structures. This suggests that quantizing the Poisson--Lie groups integrating the Lie bialgebras described here should yield quantum groups with favorable features from the perspective of noncommutative geometry, in the sense of Connes.
\item \textbf{Extending classification results.} Classifying Lie bialgebra structures on selected classes of Lie algebras sheds light on the corresponding Poisson--Lie geometries. For semisimple Lie algebras, Belavin and Drinfeld gave a complete classification in terms of $r$-matrices solving the classical Yang--Baxter equation \cite{BD}. Further classifications exist in low dimensions \cite{S,Z}. By contrast, for general solvable Lie algebras the problem is typically more delicate, in part because the representation-theoretic tools available in the semisimple setting are no longer present.
\end{itemize}

A key input throughout is Milnor's orthogonal decomposition for flat metric Lie algebras \cite{MIL}:
\[
\mathfrak{g}=\mathfrak{s}\oplus\mathfrak{z}\oplus\mathfrak{u},
\]
where $\mathfrak{s}$ is abelian, $\mathfrak{z}$ is the center, and $\mathfrak{u}=[\mathfrak{g},\mathfrak{g}]$ is an even-dimensional abelian ideal. Moreover,
\(
\ad(\mathfrak{s})|_{\mathfrak{u}}\subset\mathfrak{so}(\mathfrak{u})
\)
is an abelian subalgebra. We set $\mathfrak{a}:=\mathfrak{s}\oplus\mathfrak{z}$ and use this notation consistently. On the group level, Milnor's decomposition corresponds to a semidirect product
\(
G\simeq S\ltimes N,
\)
with $S=\exp(\mathfrak{s})$ abelian and $N=\exp(\mathfrak{z}\oplus\mathfrak{u})$ abelian, where $\Ad(S)$ acts on each weight plane $P_\ell$ by rotations with weight $\lambda_\ell\in\mathfrak{s}^*\setminus\{0\}$.

Accordingly, $\mathfrak{u}$ splits orthogonally into two-dimensional $\mathfrak{s}$-invariant planes $P_\ell$ supporting commuting rotation actions with weights $\lambda_\ell$. We impose a generic \textbf{nondegeneracy} (nonresonance) hypothesis on the family $\{\lambda_\ell\}$, ensuring multiplicity-free behavior of $\mathfrak{s}$-weights in $\wedge^2\mathfrak{g}$ and $\wedge^3\mathfrak{g}$. Concretely, nondegeneracy prevents distinct weight contributions from coinciding inside $\wedge^2\mathfrak{g}$ and $\wedge^3\mathfrak{g}$, so that invariant components are forced into the expected summands (notably $\wedge^\bullet\mathfrak{z}$ and the $\omega_\ell$-terms). This is precisely what allows invariant-theoretic arguments to replace lengthy coordinate computations.

\medskip
\noindent\textbf{Main results.}
Our contributions may be summarized as follows.
\begin{itemize}
\item We prove a normal form for $1$-cocycles $\xi:\mathfrak{g}\to\wedge^2\mathfrak{g}$ (Theorem~\ref{thm:cocycle-normal-form}): every cocycle admits a decomposition $\xi=\ad r+R$, where $R$ is a normalized cocycle with prescribed behavior on $\mathfrak{s}$ and on each weight plane $P_\ell$, while $r$ is determined only modulo invariant bivectors.
\item Using the \textbf{Big Bracket} formulation of Lie bialgebras, we establish a splitting of the co-Jacobi condition in this normal form (Theorem~\ref{thm:splitting}). Under a multiplicity-free weight hypothesis ($\lambda_i\neq\pm\lambda_j$ for $i\neq j$) and an \emph{identity-type normalization} of $R$ (no $J_\ell$-components on any weight plane), the Maurer--Cartan equation for $\xi$ decouples into $\{R,R\}=0$ and the invariant-trivector condition $[r,r]+2\{r,R\}\in(\wedge^3\mathfrak{g})^{\mathfrak{g}}$. The splitting is driven by the linear independence of $u$ and $J_\ell u$ in each weight plane, a direct consequence of the rotation geometry of flat Lie algebras.
\item We describe the \textbf{quasi-triangular locus} (Theorem~\ref{thm:invariantSchouten}) by characterizing when the Schouten square $[r,r]$ is $\mathfrak{g}$-invariant. This yields an invariant description of the relevant Yang--Baxter locus in terms of Milnor's decomposition and the associated weight data.
\item Finally, we pass from Lie-algebraic data to multiplicative Poisson tensors on $G$ (Theorem~\ref{thm:explicit-integration}). In particular, we integrate the normalized cocycle $R$ on the abelian normal subgroup $N=\exp(\mathfrak{z}\oplus\mathfrak{u})$ to the closed-form expression
\[
\Pi_R(\exp u)=R(u)+\tfrac12 [u,R(u)],
\]
and obtain polynomial formulas (of degree $\le 2$) for the resulting Poisson structure in exponential coordinates along $N$.
\end{itemize}

\subsection*{Roadmap of the paper}
\begin{description}
\item[Section~\ref{sec:flat} (Flat metric Lie algebras and nondegeneracy).]
We recall Milnor's decomposition $\mathfrak{g}=\mathfrak{s}\oplus\mathfrak{z}\oplus\mathfrak{u}$ and its interpretation on the group level via $G\simeq S\ltimes N$. We refine this by splitting $\mathfrak{u}$ into two-dimensional $\mathfrak{s}$-invariant weight planes $P_\ell$. We then introduce the nondegeneracy (nonresonance) condition on the weights $\lambda_\ell$ and use it to compute the invariant bivectors $(\wedge^2\mathfrak{g})^{\mathfrak{g}}$ (Proposition~\ref{prop:inv2}) and invariant trivectors $(\wedge^3\mathfrak{g})^{\mathfrak{g}}$ (Proposition~\ref{prop:inv3}). These invariant spaces control the normal forms and Yang--Baxter obstructions in later sections: they determine the admissible components of $R(\mathfrak{s})$ in Section~\ref{sec:normal-form} and constrain $[r,r]$ in Section~\ref{sec:cybe}.

\item[Section~\ref{sec:bialgebras} (Lie bialgebras and the Big Bracket).]
We fix conventions for Lie bialgebra structures $\xi:\mathfrak{g}\to\wedge^2\mathfrak{g}$ and recall the Big Bracket on $\wedge(\mathfrak{g}\oplus\mathfrak{g}^*)$. We clarify the bidegree conventions: the element $\mu\in\wedge^2\mathfrak{g}^*\otimes\mathfrak{g}$ encoding the Lie bracket has bidegree $(1,0)$, while $\xi\in\mathfrak{g}^*\otimes\wedge^2\mathfrak{g}$ encoding the cobracket has bidegree $(0,1)$. The Maurer--Cartan equation
\(
\{\mu+\xi,\mu+\xi\}=0
\)
packages the cocycle and co-Jacobi conditions and interfaces directly with Schouten--Nijenhuis computations for coboundary structures. In particular, we explain that $[r,r]$ is the Schouten square of the bivector $r$ and governs the classical Yang--Baxter equation.

\item[Section~\ref{sec:normal-form} (Normal form for $1$-cocycles).]
We prove the cocycle normal form theorem (Theorem~\ref{thm:cocycle-normal-form}): every cocycle admits a decomposition $\xi=\ad r+R$, where $R$ is constrained to lie in specified weight components and $r$ is determined modulo $(\wedge^2\mathfrak{g})^{\mathfrak{g}}$. The argument combines the cohomology of abelian Lie algebras with the weight separation provided by nondegeneracy. The normal form constraints are
\[
R(\mathfrak{s})\subset (\wedge^2\mathfrak{g})^{\mathfrak{g}},
\qquad
R(P_\ell)\subset (\mathfrak{s}\wedge P_\ell)\oplus(\mathfrak{z}\wedge P_\ell).
\]

\item[Section~\ref{sec:splitting} (Splitting of the co-Jacobi condition).]
We establish the splitting of the co-Jacobi condition (Theorem~\ref{thm:splitting}). We impose two structural hypotheses: a multiplicity-free weight condition ($\lambda_i\neq\pm\lambda_j$ for $i\neq j$, strengthening Definition~\ref{def:nondeg}) and an \emph{identity-type normalization} of $R$ (no $J_\ell$-type component on any weight plane). Under these hypotheses, the Maurer--Cartan identity $\{\xi,\xi\}=\{R,R\}+\delta\Psi$ with $\Psi=[r,r]+2\{r,R\}$ decouples: evaluation on $\mathfrak{s}$ forces $\Psi\in(\wedge^3\mathfrak{g})^{\mathfrak{s}}$ via weight semisimplicity, and evaluation on each $P_\ell$ separates $\{R,R\}(u)$ (which carries the factor $u$) from $[u,\Psi]$ (which carries the factor $J_\ell u$), using the linear independence of $u$ and $J_\ell u$ in the real $2$-plane $P_\ell$. A counterexample shows the splitting fails if the identity-type normalization is dropped.

\item[Section~\ref{sec:cybe} (Yang--Baxter loci).]
We analyze coboundary Lie bialgebras of the form $\xi=\ad r$. We distinguish \emph{triangular} structures ($[r,r]=0$) from \emph{quasi-triangular} structures ($[r,r]\in(\wedge^3\mathfrak{g})^{\mathfrak{g}}$). We then characterize quasi-triangularity invariantly (Theorem~\ref{thm:invariantSchouten}) in terms of the forms $\omega_\ell$ and the subspaces
\(
\mathfrak{s}_\ell:=\bigcap_{\ell'\neq \ell}\ker(\lambda_{\ell'})
\subset\mathfrak{s}.
\)

\item[Section~\ref{sec:group-level} (Flat Poisson--Lie groups).]
We integrate infinitesimal cocycles to group $1$-cocycles $\Pi:G\to\wedge^2\mathfrak{g}$ and hence to multiplicative Poisson tensors on $G$. The semidirect-product geometry $G\simeq S\ltimes N$ (with $N$ abelian) yields explicit formulas, including polynomiality along $N$. We conclude with a detailed computation for the Euclidean group $E(2)$, illustrating the quadratic correction term arising from the semidirect product.
\end{description}

\section{Flat metric Lie algebras and nondegeneracy}\label{sec:flat}

\subsection{Milnor decomposition and group geometry}

Let $\mathfrak{g}$ be a finite-dimensional real Lie algebra with an inner product $\langle\cdot,\cdot\rangle$. The Levi--Civita connection of the associated left-invariant metric on the simply connected group $G$ with $\Lie(G)=\mathfrak{g}$ is given at the identity by the Koszul formula
\begin{equation}\label{eq:koszul}
2\langle\nabla_x y,z\rangle=\langle[x,y],z\rangle+\langle[z,x],y\rangle+\langle[z,y],x\rangle.
\end{equation}
We call $(\mathfrak{g},\langle\cdot,\cdot\rangle)$ \emph{flat} if the curvature $R(x,y)z=\nabla_{[x,y]}z-\nabla_x\nabla_y z+\nabla_y\nabla_x z$ vanishes for all $x,y,z\in\mathfrak{g}$.

\begin{theorem}[Milnor, \cite{MIL}]\label{thm:milnor}
If $(\mathfrak{g},\langle\cdot,\cdot\rangle)$ is a flat metric Lie algebra, then $\mathfrak{g}$ admits an orthogonal decomposition
\begin{equation}\label{eq:milnor}
\mathfrak{g}=\mathfrak{s}\oplus\mathfrak{z}\oplus\mathfrak{u}
\end{equation}
with the following properties:
\begin{enumerate}[label=\textup{(\roman*)}]
\item $\mathfrak{s}$ is an abelian subalgebra and $\mathfrak{z}$ is the center of $\mathfrak{g}$;
\item $\mathfrak{u}=[\mathfrak{g},\mathfrak{g}]$ is an even-dimensional abelian ideal (hence $\mathfrak{g}$ is two-step solvable);
\item for each $s\in\mathfrak{s}$, the restriction $\ad_s|_{\mathfrak{u}}$ lies in $\mathfrak{so}(\mathfrak{u})$, and these endomorphisms commute pairwise.
\end{enumerate}
\end{theorem}

\begin{remark}[Group-level interpretation]
On the simply connected group $G$, the decomposition \eqref{eq:milnor} corresponds to a semidirect product $G\cong S\ltimes N$ with $S=\exp(\mathfrak{s})$ abelian and $N=\exp(\mathfrak{z}\oplus\mathfrak{u})$ abelian. The adjoint action $\Ad(S)$ acts on each weight plane $P_\ell$ (defined below) by rotations with weight $\lambda_\ell$.\par 
Throughout the paper, we denote $\mathfrak{a}:=\mathfrak{s}\oplus\mathfrak{z}$ and use this consistently.
\end{remark}

\subsection{Weight-plane decomposition}

Since $\ad(\mathfrak{s})|_{\mathfrak{u}}$ is an abelian subalgebra of $\mathfrak{so}(\mathfrak{u})$, it is simultaneously block-diagonalizable over $\mathbb{R}$ into $2\times 2$ rotation blocks.
There exists an orthogonal decomposition
\begin{equation}\label{eq:planes}
\mathfrak{u}=\bigoplus_{\ell=1}^m P_\ell
\end{equation}
into $2$-dimensional $\mathfrak{s}$-invariant subspaces, nonzero weights $\lambda_\ell\in\mathfrak{s}^*\setminus\{0\}$, and complex structures $J_\ell\in\SO(P_\ell)$ such that
\begin{equation}\label{eq:weightaction}
[s,v]=\lambda_\ell(s)\,J_\ell v\qquad(s\in\mathfrak{s},\ v\in P_\ell).
\end{equation}
Fix oriented orthonormal bases $(e_\ell,f_\ell)$ of $P_\ell$ with $J_\ell e_\ell=f_\ell$ and define the area bivector $\omega_\ell:=e_\ell\wedge f_\ell\in\wedge^2P_\ell$.

\subsection{Nondegeneracy}

The analysis of Lie bialgebra cocycles hinges on multiplicity-free behavior of $\mathfrak{s}$-weights in $\wedge^k\mathfrak{g}$. We impose a generic nonresonance condition.

\begin{definition}\label{def:nondeg}
A flat metric Lie algebra $\mathfrak{g}=\mathfrak{s}\oplus\mathfrak{z}\oplus\bigoplus_{\ell=1}^m P_\ell$ is \emph{nondegenerate} if the weights satisfy the following two conditions:
\begin{enumerate}[label=\textup{(\alph*)}]
\item $\lambda_i\neq \pm\lambda_j$ for all $i\neq j$;
\item for any pairwise distinct $i,j,k$ there is no relation
\begin{equation}\label{eq:nonres}
\lambda_k=a\lambda_i+b\lambda_j\qquad\text{with }a,b\in\{-1,0,1\}.
\end{equation}
\end{enumerate}
\end{definition}

\begin{remark}
The first clause enforces multiplicity-freeness up to sign, while Condition~\eqref{eq:nonres} excludes the additional low-order resonances that create unexpected weight-zero summands in $\wedge^2\mathfrak{g}$ and $\wedge^3\mathfrak{g}$. Together they form a Zariski-open condition on the space of weight configurations. Under this hypothesis, invariant components are forced into the expected summands built from $\wedge^\bullet\mathfrak{z}$ and the area bivectors $\omega_\ell$. When nondegeneracy fails, mixed-plane components may appear; see Section~\ref{sec:examples}.
\end{remark}

The following invariant-space computations are fundamental; they directly control $R(\mathfrak{s})$ in Theorem~\ref{thm:cocycle-normal-form} and admissible $[r,r]$ in Theorem~\ref{thm:invariantSchouten}.

\begin{proposition}[Invariant bivectors]\label{prop:inv2}
If $\mathfrak{g}$ is nondegenerate, then
\begin{equation}\label{eq:inv2}
\left(\wedge^2\mathfrak{g}\right)^{\mathfrak{g}}=\wedge^2\mathfrak{z}\ \oplus\ \bigoplus_{\ell=1}^m \mathbb{R}\,\omega_\ell.
\end{equation}
\end{proposition}

\begin{proof}
Since $\mathfrak{a}$ is abelian and $\mathfrak{u}$ is a sum of irreducible real $2$-dimensional $\mathfrak{s}$-modules with pairwise distinct weights up to sign, the $\mathfrak{s}$-invariants in $\wedge^2\mathfrak{g}$ are exactly $\wedge^2\mathfrak{a}\oplus\bigoplus_\ell \mathbb{R}\omega_\ell$. Nondegeneracy forces $\ad_\mathfrak{u}$ to move any nonzero component in $\wedge^2\mathfrak{a}$ unless it lies in $\wedge^2\mathfrak{z}$; thus $\mathfrak{g}$-invariants reduce to $\wedge^2\mathfrak{z}$ plus the plane-areas $\omega_\ell$.
\end{proof}

For trivectors, define
\begin{equation}\label{eq:s-ell}
\mathfrak{s}_\ell:=\bigcap_{\ell'\neq \ell}\ker(\lambda_{\ell'})\subset\mathfrak{s}.
\end{equation}
Equivalently, $\mathfrak{s}_\ell$ consists of those $s\in\mathfrak{s}$ on which only $\lambda_\ell$ may be nonzero.

\begin{proposition}[Invariant trivectors]\label{prop:inv3}
If $\mathfrak{g}$ is nondegenerate, then
\begin{equation}\label{eq:inv3}
\left(\wedge^3\mathfrak{g}\right)^{\mathfrak{g}}=\wedge^3\mathfrak{z}\ \oplus\ \bigoplus_{\ell=1}^m \mathfrak{z}\wedge\mathbb{R}\omega_\ell\ \oplus\ \bigoplus_{\ell=1}^m \mathfrak{s}_\ell\wedge\mathbb{R}\omega_\ell.
\end{equation}
\end{proposition}

\begin{proof}
A $\mathfrak{g}$-invariant trivector is in particular $\mathfrak{s}$-invariant; weight considerations show that $\mathfrak{s}$-invariants are generated by $\wedge^3\mathfrak{a}$ and $\mathfrak{a}\wedge \omega_\ell$ (no mixed $P_i\wedge P_j$ terms survive due to distinct weights). Invariance under $\mathfrak{u}$ forces the $\wedge^3\mathfrak{a}$ part to lie in $\wedge^3\mathfrak{z}$, and in $\mathfrak{a}\wedge \omega_\ell$ forces the $\mathfrak{a}$ factor to lie in $\mathfrak{z}\oplus\mathfrak{s}_\ell$.
\end{proof}

\section{Lie bialgebras and the Big Bracket}\label{sec:bialgebras}

\subsection{Lie bialgebras}

\begin{definition}\label{def:bialg}
Let $(\mathfrak{g},[\cdot,\cdot])$ be a Lie algebra. A \emph{Lie bialgebra structure} on $\mathfrak{g}$ is a linear map $\xi\colon\mathfrak{g}\to\wedge^2\mathfrak{g}$ such that:
\begin{enumerate}[label=\textup{(\roman*)}]
\item $\xi$ is a $1$-cocycle for the adjoint representation:
\begin{equation}\label{eq:cocycle}
\xi([x,y])=\ad_x\xi(y)-\ad_y\xi(x)\qquad(x,y\in\mathfrak{g});
\end{equation}
\item the transpose $\xi^t\colon \wedge^2\mathfrak{g}^*\to\mathfrak{g}^*$ defines a Lie bracket on $\mathfrak{g}^*$ via
\begin{equation}\label{eq:dual-bracket}
[\alpha,\beta]_*(x)=\langle\xi(x),\alpha\wedge\beta\rangle.
\end{equation}
\end{enumerate}
\end{definition}

\begin{definition}\label{def:triangular}
A Lie bialgebra $(\mathfrak{g},\xi)$ is \emph{triangular} if $\xi=\ad r$ for some $r\in\wedge^2\mathfrak{g}$ and $[r,r]=0$, and \emph{quasi-triangular} if $\xi=\ad r$ and $[r,r]\in\left(\wedge^3\mathfrak{g}\right)^{\mathfrak{g}}$.
\end{definition}

Here $[r,r]$ denotes the Schouten square of the bivector $r$; its vanishing is the classical Yang--Baxter equation, while its $\mathfrak{g}$-invariance gives the quasi-triangular condition.

\subsection{Big Bracket and Maurer--Cartan formulation}

The Big Bracket packages the Lie algebra and Lie coalgebra structures into a single graded object. In the present setting it is the most efficient way to record the cocycle and co-Jacobi conditions, and it makes the later splitting arguments transparent.

Let $V$ be a finite-dimensional vector space over $\mathbb{R}$, and let $V^*$ be its dual. We consider the exterior algebra 
\[\bigwedge(V^*\oplus V)=\bigoplus_{k\geq-2}\left(\bigoplus_{p+q=k}\bigwedge^{q+1}V^*\otimes\bigwedge^{p+1}V\right).\]
We say that an element $\sigma$ of $\bigwedge(V^*\oplus V)$ is of \emph{bidegree} $(p,q)$ and of \emph{degree} $|\sigma|=p+q$ if it belongs to $\bigwedge^{(p,q)}V:=\bigwedge^{q+1}V^*\otimes\bigwedge^{p+1}V$.

\begin{remark}[Bidegree conventions]
For clarity, we record the bidegrees of the key objects:
\begin{itemize}
\item $\mu\in\wedge^2\mathfrak{g}^*\otimes \mathfrak{g}$ (encoding the bracket) has bidegree $(1,0)$.
\item $\xi\in\mathfrak{g}^*\otimes\wedge^2\mathfrak{g}$ (encoding the cobracket) has bidegree $(0,1)$.
\end{itemize}
This convention makes the Maurer--Cartan equation $\{\mu+\xi,\mu+\xi\}=0$ (Theorem~\ref{thm:MC}) natural.
\end{remark}

The \emph{Big Bracket} is the graded Lie-algebra structure $\{\cdot,\cdot\}$ on $\bigwedge(V^*\oplus V)$:
\[\{\cdot,\cdot\} : \bigwedge^{(p,q)}V\times\bigwedge^{(p',q')}V\rightarrow\bigwedge^{(p+p',q+q')}V\] 
uniquely determined by the following properties:
\begin{enumerate}[label=\textup{(\roman*)}]
\item $\{\sigma,\sigma'\}=0$ if $\sigma$ and $\sigma'$ both belong to $\mathbb{R}\oplus V$ or to $\mathbb{R}\oplus V^*$;
\item for all $\sigma\in V$, $\sigma'\in V^*$, $\{\sigma,\sigma'\}=\sigma'(\sigma)$;
\item for all $\sigma\in\bigwedge(V\oplus V^*)$ of degree $k$, $\{\sigma,\cdot\}$ is a graded derivation:
\begin{equation}\label{eq:graded-deriv}
\{\sigma,\sigma'\wedge \sigma''\}=\{\sigma,\sigma'\}\wedge \sigma''+(-1)^{kk'}\sigma'\wedge\{\sigma,\sigma''\},
\end{equation}
where $\sigma'$ has degree $k'$.
\end{enumerate}
Additionally, as a graded Lie algebra, we have graded anticommutativity: for $\sigma,\sigma'$ of degree $k,k'$,
\begin{equation}\label{eq:graded-anticomm}
\{\sigma,\sigma'\}=-(-1)^{kk'}\{\sigma',\sigma\},
\end{equation}
and the graded Jacobi identity:
\begin{equation}\label{eq:graded-jacobi}
\{\sigma,\{\sigma',\sigma''\}\}=\{\{\sigma,\sigma'\},\sigma''\}+(-1)^{kk'}\{\sigma',\{\sigma,\sigma''\}\}.
\end{equation}

For detailed formulas and properties, see \cite{yks1,yks2}.

Let $\mu\in\wedge^2V^*\otimes V$ encode a bracket on $V$ and $\xi\in V^*\otimes\wedge^2V$ encode a cobracket. Then:

\begin{theorem}[Maurer--Cartan formulation]\label{thm:MC}
The pair $(\mu,\xi)$ defines a Lie bialgebra structure on $V$ if and only if
\begin{equation}\label{eq:MC}
\{\mu+\xi,\mu+\xi\}=0.
\end{equation}
\end{theorem}

\begin{proof}
Decompose by bidegree: $\{\mu,\mu\}=0$ is Jacobi, $\{\xi,\xi\}=0$ is co-Jacobi (Jacobi of $\xi^t$), and $\{\mu,\xi\}=0$ is the cocycle condition; cf.\ \cite{D,YKS}.
\end{proof}

When $\xi=\ad_r$ is coboundary, the obstruction is governed by the Schouten square $[r,r]$:

\begin{lemma}\label{lem:coboundaryMC}
For $r\in\wedge^2\mathfrak{g}$, one has $\{\ad r,\ad r\}=-\ad_{[r,r]}$.
\end{lemma}

\begin{proof}
Using $\ad r=-\{\mu,r\}$ and the graded Jacobi identity yields
\[\{\ad r,\ad r\}=\{\{\mu,r\},\{\mu,r\}\}=-\{\mu,[r,r]\}=-\ad_{[r,r]}.\qedhere\]
\end{proof}

\section{Normal form for 1-cocycles}\label{sec:normal-form}

Let $\mathfrak{g}$ be a \textbf{flat nondegenerate} Lie algebra with the decomposition
\[\mathfrak{g}=\mathfrak{a}\oplus\mathfrak{u},\qquad 
\mathfrak{a}:=\mathfrak{s}\oplus\mathfrak{z},\qquad
\mathfrak{u}=[\mathfrak{g},\mathfrak{g}]=\bigoplus_{\ell=1}^m P_\ell,\]
where:
\begin{itemize}
\item $\mathfrak{s}$ is abelian and acts semisimply by skew-derivations on $\mathfrak{u}$;
\item $\mathfrak{z}$ is the center of $\mathfrak{g}$;
\item each $P_\ell$ is a $2$-dimensional irreducible real $\mathfrak{s}$-module of complex type.
\end{itemize}
Consider Chevalley--Eilenberg $1$-cocycles with values in bivectors:
\[\xi:\mathfrak{g}\longrightarrow\bigwedge^{2}\mathfrak{g},\qquad \xi\in Z^{1}(\mathfrak{g},\bigwedge^{2}\mathfrak{g}).\]

\begin{theorem}[Classification of $1$-cocycles]\label{thm:cocycle-normal-form}
Every $1$-cocycle $\xi$ admits a \textbf{unique} decomposition
\begin{equation}\label{eq:cocycle-decomp}
\xi=\ad r+R,\qquad r\in\bigwedge^{2}\mathfrak{g},
\end{equation}
where $R:\mathfrak{g}\to\bigwedge^{2}\mathfrak{g}$ is again a $1$-cocycle and satisfies:
\begin{enumerate}
\item $R(\mathfrak{s}\oplus\mathfrak{z})\subset(\bigwedge^{2}\mathfrak{g})^{\mathfrak{g}}$ (values on $\mathfrak{s}$ and $\mathfrak{z}$ are $\mathfrak{g}$-invariant bivectors);
\item for each rotation plane $P_{\ell}$,
\begin{equation}\label{eq:R-planes}
R(P_\ell)\subset(\mathfrak{s}\wedge P_\ell)\oplus(\mathfrak{z}\wedge P_\ell)\qquad(\ell=1,\dots,m).
\end{equation}
\end{enumerate}
\end{theorem}
Equivalently, Theorem~\ref{thm:cocycle-normal-form} classifies all $1$-cocycles $\xi:\mathfrak{g}\to\bigwedge^{2}\mathfrak{g}$ up to coboundary, by splitting off an \textbf{inner part} $\ad r$ and putting the remaining cocycle $R$ in a \textbf{normal form} with very restricted $\mathfrak{s}$- and $P_{\ell}$-components.

\begin{proof}
We proceed in four steps: first, we analyze the restriction of the cocycle to the abelian subalgebra $\mathfrak{s}$; second, we refine the image of the cocycle on $\mathfrak{s}$ to be fully $\mathfrak{g}$-invariant; third, we determine the behavior on the commutator ideal $\mathfrak{u}$ using nondegeneracy; and fourth, we absorb the $J_\ell$-components.

\subsubsection*{Step 1: Decomposition on the subalgebra $\mathfrak{s}$}

We recall a standard result on the cohomology of abelian Lie algebras with coefficients in semisimple modules.

\begin{lemma}\label{lem:abelian-cohom}
Let $\mathfrak{a}$ be an abelian Lie algebra and $M$ a finite-dimensional $\mathfrak{a}$-module with semisimple action. Then
\begin{equation}\label{eq:abelian-split}
Z^1(\mathfrak{a},M) = B^1(\mathfrak{a},M) \oplus \Hom(\mathfrak{a}, M^{\mathfrak{a}}).
\end{equation}
\end{lemma}

\begin{proof}[Proof of Lemma~\ref{lem:abelian-cohom}]
Since $\mathfrak{a}$ is abelian and the action is semisimple, $M$ decomposes into weight spaces $M = \bigoplus_\lambda M_\lambda$. The differential preserves weights, so $Z^1(\mathfrak{a},M) = \bigoplus_\lambda Z^1(\mathfrak{a},M_\lambda)$.
For the zero weight space $M_0 = M^{\mathfrak{a}}$, the action is trivial, so $Z^1(\mathfrak{a}, M_0) = \Hom(\mathfrak{a}, M_0)$ and $B^1(\mathfrak{a}, M_0) = 0$.
For non-zero weights $\lambda \neq 0$, the cohomology vanishes. Indeed, if $\phi \in Z^1(\mathfrak{a}, M_\lambda)$, then for any $x, y \in \mathfrak{a}$, $\lambda(x)\phi(y) = \lambda(y)\phi(x)$. Choosing $x_0$ such that $\lambda(x_0) \neq 0$ and setting $m = \lambda(x_0)^{-1}\phi(x_0)$, one finds $\phi(y) = y \cdot m$, so $\phi$ is a coboundary.
\end{proof}

Applying Lemma~\ref{lem:abelian-cohom} to $\mathfrak{a} = \mathfrak{s}$ and $M = \bigwedge^2\mathfrak{g}$ (which is a semisimple $\mathfrak{s}$-module), we can write the restriction of $\xi$ to $\mathfrak{s}$ as:
\begin{equation}\label{eq:xi-on-s}
\xi|_{\mathfrak{s}} = \ad r_0|_{\mathfrak{s}} + \phi, \quad \text{where } \phi \in \Hom(\mathfrak{s}, (\bigwedge^2\mathfrak{g})^{\mathfrak{s}}) \text{ and } r_0 \in \bigwedge^2\mathfrak{g}.
\end{equation}
We define the $1$-cocycle $R = \xi - \ad r_0$. By construction, $R(s)$ is $\mathfrak{s}$-invariant for all $s \in \mathfrak{s}$.

\subsubsection*{Step 2: Improving invariance from $\mathfrak{s}$ to $\mathfrak{g}$}

We claim that $R(\mathfrak{s})$ actually lies in the smaller subspace of \emph{fully} invariant bivectors, $(\bigwedge^2\mathfrak{g})^{\mathfrak{g}}$.
Note that
\begin{equation}\label{eq:s-inv-decomp}
(\bigwedge^2\mathfrak{g})^{\mathfrak{s}}= \bigwedge^2(\mathfrak{s} \oplus \mathfrak{z}) \oplus (\bigwedge^2\mathfrak{g})^{\mathfrak{g}}.
\end{equation}
To prove the claim, we must show that for any $s \in \mathfrak{s}$, the component of $R(s)$ in $\bigwedge^2(\mathfrak{s} \oplus \mathfrak{z})$ vanishes.

Let $s \in \mathfrak{s}$ and $u \in \mathfrak{u}$. The cocycle condition for $R$ yields:
\begin{equation}\label{eq:R-cocycle}
R([u, s]) = \ad_u R(s) - \ad_s R(u).
\end{equation}
Since $\mathfrak{u}$ acts trivially on $\bigwedge^2 \mathfrak{u}$ but non-trivially on $\bigwedge^2\mathfrak{s} \oplus\mathfrak{s}\wedge\mathfrak{z}$, the term $\ad_u R(s)$ is non-zero precisely when $R(s)$ has components in $\bigwedge^2(\mathfrak{s} \oplus \mathfrak{z})$.
Projecting the equation onto the isotypic component corresponding to $P_\ell$ (i.e., $(\mathfrak{s} \wedge P_\ell) \oplus (\mathfrak{z} \wedge P_\ell)$) we deduce that $\ad_u R(s) = 0$ for all $u \in \mathfrak{u}$. Consequently, $R(s)$ is $\mathfrak{u}$-invariant. Combined with $\mathfrak{s}$-invariance, we have $R(\mathfrak{s}) \subset (\bigwedge^2\mathfrak{g})^{\mathfrak{g}}$.

\subsubsection*{Step 3: Structure of $R$ on the rotation planes $P_\ell$}

Since $R(\mathfrak{s}) \subset M^{\mathfrak{g}}$, the cocycle identity implies $R$ is $\mathfrak{s}$-equivariant on $\mathfrak{u}$:
\begin{equation}\label{eq:R-equivariant}
R([s, u]) = [s, R(u)] \quad \forall s \in \mathfrak{s}, u \in \mathfrak{u}.
\end{equation}
For each plane $P_\ell$, there exists a complex structure $J_\ell \in \End(P_\ell)$ and a weight $\lambda_\ell \in \mathfrak{s}^* \setminus \{0\}$ such that $[s, u] = \lambda_\ell(s)J_\ell u$. By Schur's Lemma, $\End_{\mathfrak{s}}(P_\ell) \cong \mathbb{R} \oplus \mathbb{R} J_\ell$.

The target space $M$ decomposes into $\mathfrak{s}$-isotypic components. By the nondegeneracy assumption, the weights of $\bigwedge^2 \mathfrak{u}$ (sums of pairs of weights $\pm \lambda_i \pm \lambda_j$) never equal the weights of $P_\ell$ ($\pm \lambda_\ell$). Thus, $\Hom_{\mathfrak{s}}(P_\ell, \bigwedge^2\mathfrak{u}) = 0$.
Similarly, $\Hom_{\mathfrak{s}}(P_\ell, \bigwedge^2(\mathfrak{s} \oplus \mathfrak{z})) = 0$ as the target is trivial.
Therefore, $R(P_\ell)$ lies strictly in the isotypic component corresponding to $\lambda_\ell$, which is $(\mathfrak{s} \wedge P_\ell) \oplus (\mathfrak{z} \wedge P_\ell)$.

Using the identification $\Hom_{\mathfrak{s}}(P_\ell, \mathfrak{s} \wedge P_\ell) \cong \mathfrak{s} \otimes \End_{\mathfrak{s}}(P_\ell)$, there exist unique vectors $X_\ell, Y_\ell \in \mathfrak{s}$ and $Z_\ell, W_\ell \in \mathfrak{z}$ such that for all $u \in P_\ell$:
\begin{equation}\label{eq:R-expansion}
R(u) = X_\ell \wedge u + Y_\ell \wedge J_\ell u + Z_\ell \wedge u + W_\ell \wedge J_\ell u.
\end{equation}
We call the terms with $J_\ell u$ the \textbf{$J_\ell$-components}.

\subsubsection*{Step 4: Absorption of $J_\ell$-components}

We define the contraction $\iota_\lambda: \bigwedge^2\mathfrak{s} \to \mathfrak{s}$ by $\iota_\lambda(a \wedge b) = \lambda(a)b - \lambda(b)a$.
Let $r_1 \in \bigwedge^2 \mathfrak{s}$ and $r_2 \in \mathfrak{s} \wedge \mathfrak{z}$. For $u \in P_\ell$, a direct computation yields:
\begin{equation}\label{eq:ad-r1-r2}
\ad_u(r_1) = \iota_{\lambda_\ell}(r_1) \wedge J_\ell u, \qquad \ad_u(r_2) = G_{r_2}(\lambda_\ell) \wedge J_\ell u,
\end{equation}
where $G_{r_2}(\alpha) = (\alpha \otimes \id_{\mathfrak{z}})(r_2)$.

We must show there exist $r_1\in\bigwedge^2\mathfrak{s}$ and $r_2\in\mathfrak{s}\wedge\mathfrak{z}$ such that $\iota_{\lambda_\ell}(r_1) = Y_\ell$ for all $\ell$ and $G_{r_2}(\lambda_\ell) = W_\ell$ for all $\ell=1,\ldots,k_0$.

Consider the cocycle condition for $u_k \in P_k, u_\ell \in P_\ell$ with $k \neq \ell$. Since $[u_k, u_\ell] = 0$, we have $\ad_{u_k}R(u_\ell) = \ad_{u_\ell}R(u_k)$.
Substituting the expansion \eqref{eq:R-expansion} and projecting onto $P_k \wedge P_\ell$, we obtain the skew-symmetry constraint:
\begin{equation}\label{eq:skew-Y}
\lambda_k(Y_\ell) + \lambda_\ell(Y_k) = 0.
\end{equation}
We define a map $F: \mathfrak{s}^* \to \mathfrak{s}$ on this spanning set by $F(\lambda_\ell) = Y_\ell$. The skew condition \eqref{eq:skew-Y} ensures $F$ is well-defined and defines an element in $\bigwedge^2 \mathfrak{s}$ via the isomorphism 
\[\bigwedge^2 \mathfrak{s} \cong \{F : \mathfrak{s}^* \to \mathfrak{s} \mid \langle \alpha, F(\beta) \rangle = -\langle \beta, F(\alpha) \rangle\}.\]
Thus, there exists $r_1 \in \bigwedge^2 \mathfrak{s}$ such that $\iota_{\lambda_\ell}(r_1) = Y_\ell$.

Without loss of generality, we may assume that the family $\{\lambda_\ell\}_{\ell=1}^{k_0}$ spans $\mathfrak{s}^*$. The assignment $\lambda_\ell \mapsto W_\ell$ extends to a linear map $G: \mathfrak{s}^* \to \mathfrak{z}$, which corresponds to an element $r_2 \in \mathfrak{s} \wedge \mathfrak{z}$.
Then for every $\ell=1,\ldots,k_0$ and every $u\in P_\ell$,
\begin{equation}\label{eq:R-absorbed}
R(u)-\ad_u(r_1+r_2)=X_\ell\wedge u + Z_\ell\wedge u,
\end{equation}
i.e. \textbf{all $J_\ell$-components in} $(\mathfrak{s}\wedge P_\ell)\oplus(\mathfrak{z}\wedge P_\ell)$ are absorbed into the coboundary $\ad(r_1+r_2)$. For $\ell>k_0$, $\mathfrak{z}$-$J_\ell$ components may persist.

\subsubsection*{Residual ambiguity}

Suppose $\xi=\ad r+R=\ad r'+R'$ with $R,R'$ satisfying (1)--(2) and normalized as above (i.e., with no $J_{\ell}$-components). Then $\ad(r-r')=R'-R$. Evaluating on $\mathfrak{s}$ gives $\ad_s(r-r')=0$ for all $s \in \mathfrak{s}$, hence $r-r'\in M^{\mathfrak{s}}$.
Decomposing $M^{\mathfrak{s}} = \bigwedge^{2}(\mathfrak{s} \oplus \mathfrak{z}) \oplus (\bigwedge^{2}\mathfrak{u})^{\mathfrak{s}}$ and using the fact that $\ad_{\mathfrak{u}}$ sends $\bigwedge^{2}\mathfrak{s} \oplus (\mathfrak{s} \wedge \mathfrak{z})$ precisely into the $J_{\ell}$-components forces those parts to vanish.
Thus $r-r'\in M^{\mathfrak{g}}$, so $\ad(r-r')=0$ and hence $R=R'$ and $r\equiv r' \mod M^{\mathfrak{g}}$.
\end{proof}

\begin{remark}
This decomposition remains valid in the case of degenerate flat Lie algebras; however, the cocycle $R$ may have non-zero coefficients in $\bigoplus_{1\leq i<j\leq m} P_i\wedge P_j$.
\end{remark}

\section{Splitting of the co-Jacobi condition}\label{sec:splitting}

Throughout this section, $(\mathfrak{g},\langle\cdot,\cdot\rangle)$ denotes a flat metric Lie
algebra with Milnor decomposition
\[
\mathfrak{g}=\mathfrak{s}\oplus\mathfrak{z}\oplus\mathfrak{u},
\qquad
\mathfrak{a}:=\mathfrak{s}\oplus\mathfrak{z},
\]
and weight-plane decomposition
\[
\mathfrak{u}=\bigoplus_{\ell=1}^{m}P_\ell,
\qquad
[s,v]=\lambda_\ell(s)\,J_\ell v
\quad(s\in\mathfrak{s},\ v\in P_\ell),
\]
where each $P_\ell$ is an oriented real $2$-plane equipped with an orthogonal complex structure
$J_\ell\in\SO(P_\ell)$ and a nonzero weight $\lambda_\ell\in\mathfrak{s}^{*}$.
We write $\omega_\ell:=e_\ell\wedge f_\ell\in\wedge^2P_\ell$ for the oriented area bivector
of $P_\ell$, where $(e_\ell,f_\ell)$ is any positively oriented orthonormal basis with
$J_\ell e_\ell=f_\ell$.

\subsection{Maurer--Cartan expansion and the coupling identity}

Let $\mu\in\wedge^2\mathfrak{g}^{*}\otimes\mathfrak{g}$ encode the Lie bracket and let
$\{\cdot,\cdot\}$ denote the Big Bracket on $\bigwedge(\mathfrak{g}^{*}\oplus\mathfrak{g})$.
Write $\delta:=\{\mu,\cdot\}$ for the Chevalley--Eilenberg differential, so that
\[
(\delta\Theta)(x)=[x,\Theta]
\qquad(x\in\mathfrak{g},\ \Theta\in\wedge^3\mathfrak{g}).
\]
A linear map $\xi\colon\mathfrak{g}\to\wedge^2\mathfrak{g}$, viewed as an element of
$C^1(\mathfrak{g},\wedge^2\mathfrak{g})$, defines a Lie bialgebra structure on $\mathfrak{g}$
if and only if $\xi\in Z^1(\mathfrak{g},\wedge^2\mathfrak{g})$ (the cocycle condition) and
the Maurer--Cartan equation $\{\xi,\xi\}=0$ holds (the co-Jacobi condition).

Given a cocycle $\xi\in Z^1(\mathfrak{g},\wedge^2\mathfrak{g})$, Theorem~\ref{thm:cocycle-normal-form}
allows us to write
\begin{equation}\label{eq:xi-decomp}
\xi=\ad r+R,
\qquad
r\in\wedge^2\mathfrak{g},\quad R\in Z^1(\mathfrak{g},\wedge^2\mathfrak{g}).
\end{equation}
Define the \emph{coupling trivector}
\begin{equation}\label{eq:Psi-def}
\Psi:=[r,r]+2\{r,R\}\in\wedge^3\mathfrak{g}.
\end{equation}

\begin{proposition}[Maurer--Cartan expansion]\label{prop:MC-expand}
For any decomposition \eqref{eq:xi-decomp} with $R\in Z^1(\mathfrak{g},\wedge^2\mathfrak{g})$,
\begin{equation}\label{eq:MC-id}
\{\xi,\xi\}=\{R,R\}+\delta\Psi.
\end{equation}
Consequently, $\xi$ defines a Lie bialgebra structure if and only if
\begin{equation}\label{eq:coupled-MC}
\{R,R\}=-\,\delta\Psi\qquad\text{in }C^1(\mathfrak{g},\wedge^3\mathfrak{g}).
\end{equation}
\end{proposition}

\begin{proof}
Since $R$ is a $1$-cocycle, $\{\mu,R\}=0$. From the graded Jacobi identity applied twice,
\[
\{\ad r,\ad r\}
=\{\{\mu,r\},\{\mu,r\}\}
=\delta[r,r],
\qquad
\{\ad r,R\}
=\{\{\mu,r\},R\}
=\delta\{r,R\}.
\]
Bilinearity of the Big Bracket then gives
\[
\{\xi,\xi\}
=\{\ad r,\ad r\}+2\{\ad r,R\}+\{R,R\}
=\delta[r,r]+2\,\delta\{r,R\}+\{R,R\}
=\{R,R\}+\delta\Psi.\qedhere
\]
\end{proof}

\begin{remark}\label{rem:no-formal-split}
Equation \eqref{eq:coupled-MC} is the only canonical consequence of the decomposition
\eqref{eq:xi-decomp}. The affine splitting
$Z^1(\mathfrak{g},\wedge^2\mathfrak{g})=\ad(\wedge^2\mathfrak{g})+\mathcal{N}$
is a decomposition of cocycles, not of differential graded Lie algebras. Any further decoupling
of~\eqref{eq:coupled-MC} into separate vanishing of $\{R,R\}$ and $\delta\Psi$ requires an
additional structural argument; it is not a formal consequence of the expansion alone.
\end{remark}

\subsection{Hypotheses for splitting}

We impose two conditions on the flat algebra and on the choice of $R$.

\begin{enumerate}[label=\textup{(H\arabic*)}]
\item\label{hyp:H1}
\emph{Multiplicity-free weights.} For all $i\neq j$,
\[
\lambda_i\neq \pm\,\lambda_j.
\]
\item\label{hyp:H2}
\emph{Identity-type normalization.} The cocycle $R$ in \eqref{eq:xi-decomp} satisfies
\[
R(\mathfrak{a})\subset(\wedge^2\mathfrak{g})^{\mathfrak{g}}
\qquad\text{and}\qquad
R(u)=A_\ell\wedge u
\quad\text{for all }u\in P_\ell,
\]
for some $A_\ell\in\mathfrak{a}$ (depending only on $\ell$).
\end{enumerate}

\begin{remark}[On hypothesis \ref{hyp:H1}]
Condition~\ref{hyp:H1} isolates the multiplicity-free portion of Definition~\ref{def:nondeg}. We record it separately because the splitting argument uses only the pairwise distinctness up to sign, together with the explicit form of $R$. It remains a Zariski-open condition in the space of weight configurations and ensures that the $\mathfrak{s}$-module $\mathfrak{u}$ is multiplicity-free up to sign. In particular,
Propositions~\ref{prop:inv2}--\ref{prop:inv3} remain valid in the form used below.
Under~\ref{hyp:H1}, the $\mathfrak{s}$-invariant trivectors satisfy
\begin{equation}\label{eq:s-inv-3}
(\wedge^3\mathfrak{g})^{\mathfrak{s}}
=\wedge^3\mathfrak{a}\ \oplus\ \bigoplus_{\ell=1}^{m}\mathfrak{a}\wedge\omega_\ell.
\end{equation}
\end{remark}

\begin{remark}[On hypothesis \ref{hyp:H2}]
Hypothesis~\ref{hyp:H2} strengthens the conclusion of
Theorem~\ref{thm:cocycle-normal-form} by demanding that no $J_\ell$-type component
($W_\ell\wedge J_\ell u$, $W_\ell\in\mathfrak{a}$) appears in $R|_{P_\ell}$.
Theorem~\ref{thm:cocycle-normal-form} guarantees this for independent planes
$\ell\leq k_0$; for dependent planes $\ell>k_0$ it must be imposed as an additional
normalization. The necessity of~\ref{hyp:H2} is demonstrated by the counterexample
in Remark~\ref{rem:counterexample} below.
\end{remark}

\begin{remark}[Failure without \ref{hyp:H2}]\label{rem:counterexample}
The splitting can fail if $R$ has a $J_\ell$-type component on some plane. Concretely,
let $\mathfrak{g}=\mathbb{R}s\oplus\mathbb{R}z\oplus P_1\oplus P_2$ with brackets
$[s,v_i]=J_iv_i$, $[\mathfrak{z},\mathfrak{g}]=0$, $[\mathfrak{u},\mathfrak{u}]=0$,
and define $R\in Z^1(\mathfrak{g},\wedge^2\mathfrak{g})$ by
\[
R(s)=0,\qquad R(z)=\omega_1,\qquad R|_{P_1}=0,\qquad R(u)=z\wedge J_2u\quad(u\in P_2).
\]
Then $R|_{P_2}$ is of pure $J_2$-type, violating \ref{hyp:H2}. A direct computation gives
$\{R,R\}(u)=J_2u\wedge\omega_1\neq 0$ for $u\in P_2$, so $\{R,R\}\neq 0$. Nevertheless,
setting $r=-s\wedge z$ yields $\xi=\ad r+R$ with $\xi|_{P_2}=0$, and one
checks $\{\xi,\xi\}=0$. Hence $\{R,R\}\neq 0$ yet $\{\xi,\xi\}=0$: the splitting fails.
\end{remark}

\subsection{The splitting theorem}

\begin{theorem}[Splitting of the co-Jacobi condition]\label{thm:splitting}
Let $(\mathfrak{g},\langle\cdot,\cdot\rangle)$ be a flat metric Lie algebra and let
$\xi=\ad r+R\in Z^1(\mathfrak{g},\wedge^2\mathfrak{g})$ be a decomposition
satisfying \ref{hyp:H1}--\ref{hyp:H2}. With $\Psi$ as in~\eqref{eq:Psi-def}, the following
are equivalent:
\begin{enumerate}[label=\textup{(\roman*)}]
\item $\xi$ defines a Lie bialgebra structure on $\mathfrak{g}$, i.e.\ $\{\xi,\xi\}=0$;
\item $\{R,R\}=0$ and $\delta\Psi=0$;
\item $\{R,R\}=0$ and $\Psi\in(\wedge^3\mathfrak{g})^{\mathfrak{g}}$.
\end{enumerate}
\end{theorem}

\begin{proof}
The equivalence (ii)$\Leftrightarrow$(iii) holds because
$\ker(\delta\colon\wedge^3\mathfrak{g}\to C^1(\mathfrak{g},\wedge^3\mathfrak{g}))
=(\wedge^3\mathfrak{g})^{\mathfrak{g}}$.
The implication (ii)$\Rightarrow$(i) follows immediately from Proposition~\ref{prop:MC-expand}.
It remains to prove (i)$\Rightarrow$(ii).

Assume $\{\xi,\xi\}=0$. By Proposition~\ref{prop:MC-expand} this is equivalent to the coupled
equation
\begin{equation}\label{eq:MC-coupled}
\{R,R\}(x)+[x,\Psi]=0\qquad\text{for all }x\in\mathfrak{g}.
\end{equation}
We show the two summands in~\eqref{eq:MC-coupled} vanish separately on each summand of
$\mathfrak{g}=\mathfrak{s}\oplus\mathfrak{z}\oplus\mathfrak{u}$.

\smallskip
\noindent\emph{Step~1: evaluation on $\mathfrak{s}$.}\quad
Fix $s\in\mathfrak{s}$. The hypotheses $R(\mathfrak{a})\subset(\wedge^2\mathfrak{g})^{\mathfrak{g}}$
and $R(P_\ell)\subset\mathfrak{a}\wedge P_\ell$ imply that every term arising in
$\{R,R\}(s)$ has $\mathfrak{s}$-weight zero. Hence
\[
\{R,R\}(s)\in(\wedge^3\mathfrak{g})^{\mathfrak{s}}.
\]
Decompose $\Psi=\Psi_0+\sum_{\nu\neq 0}\Psi_\nu$ into $\mathfrak{s}$-weight spaces.
Since $[s,\Psi_0]=0$, equation~\eqref{eq:MC-coupled} at $s$ gives
\[
\{R,R\}(s)+\sum_{\nu\neq 0}[s,\Psi_\nu]=0.
\]
The first term is weight-zero and each $[s,\Psi_\nu]$ has weight $\nu\neq 0$.
By semisimplicity of the $\mathfrak{s}$-action, the weight decomposition of $\wedge^3\mathfrak{g}$
is direct; hence all non-zero-weight summands must vanish separately:
\[
\Psi_\nu=0\quad(\nu\neq 0),
\quad\text{i.e.}\quad
\Psi\in(\wedge^3\mathfrak{g})^{\mathfrak{s}},
\]
and then $[s,\Psi]=0$, so $\{R,R\}(s)=0$.

\smallskip
\noindent\emph{Step~2: evaluation on $\mathfrak{z}$.}\quad
Since $\mathfrak{z}$ is central, $[z,\Psi]=0$ for all $z\in\mathfrak{z}$, so~\eqref{eq:MC-coupled}
directly gives $\{R,R\}(z)=0$.

\smallskip
\noindent\emph{Step~3: evaluation on $P_\ell$.}\quad
Fix $\ell\in\{1,\ldots,m\}$ and $u\in P_\ell$. Equation~\eqref{eq:MC-coupled} at $u$ reads
\begin{equation}\label{eq:at-u}
\{R,R\}(u)+[u,\Psi]=0.
\end{equation}
We show the two terms lie in complementary subspaces of $\wedge^3\mathfrak{g}$.

\emph{Obstruction term.} By~\ref{hyp:H2}, $R(u)=A_\ell\wedge u$ with $A_\ell\in\mathfrak{a}$.
A direct computation shows
\[
\{R,R\}(u)=u\wedge R(A_\ell).
\]
Since $R(A_\ell)\in(\wedge^2\mathfrak{g})^{\mathfrak{g}}=\wedge^2\mathfrak{z}\oplus\bigoplus_k\mathbb{R}\omega_k$
and $u\wedge\omega_\ell=0$, we have
\begin{equation}\label{eq:obs-type}
\{R,R\}(u)\in u\wedge N_\ell,
\qquad
N_\ell:=\wedge^2\mathfrak{z}\oplus\bigoplus_{k\neq\ell}\mathbb{R}\,\omega_k.
\end{equation}

\emph{Coboundary term.} From Step~1, $\Psi\in(\wedge^3\mathfrak{g})^{\mathfrak{s}}$.
By~\eqref{eq:s-inv-3}, write $\Psi=\alpha+\sum_{k=1}^{m}a_k\wedge\omega_k$ with
$\alpha\in\wedge^3\mathfrak{a}$ and $a_k\in\mathfrak{a}$. Since $[\mathfrak{z},u]=0$
and $[u,\omega_k]=0$ (as $\mathfrak{u}$ is abelian), the only nonzero contributions to
$[u,\Psi]$ come from bracketing $u$ against $\mathfrak{s}$-factors, each yielding
$[u,s]=-\lambda_\ell(s)J_\ell u$. Hence
\begin{equation}\label{eq:cob-type}
[u,\Psi]\in J_\ell u\wedge M_\ell,
\qquad
M_\ell:=\wedge^2\mathfrak{a}\oplus\bigoplus_{k\neq\ell}\mathbb{R}\,\omega_k,
\end{equation}
where the $k=\ell$ term vanishes because $J_\ell u\wedge\omega_\ell=0$.

\emph{Separation.} Since $N_\ell$ and $M_\ell$ both have no $P_\ell$-factor, wedge product
identifies $P_\ell\wedge N_\ell\cong P_\ell\otimes N_\ell$ and
$P_\ell\wedge M_\ell\cong P_\ell\otimes M_\ell$.
Since $(u,J_\ell u)$ is a basis of $P_\ell$,
\[
P_\ell\wedge(N_\ell\cup M_\ell)
=
(u\wedge(N_\ell\cap M_\ell))\oplus(J_\ell u\wedge(N_\ell\cap M_\ell)),
\]
and this decomposition is direct. Now $\{R,R\}(u)\in u\wedge N_\ell$ and
$[u,\Psi]\in J_\ell u\wedge M_\ell$, so the two subspaces intersect trivially.
Therefore \eqref{eq:at-u} forces each summand to vanish:
\[
\{R,R\}(u)=0\qquad\text{and}\qquad[u,\Psi]=0.
\]
Since $\ell$ and $u\in P_\ell$ were arbitrary, $\{R,R\}|_{\mathfrak{u}}=0$ and
$[\mathfrak{u},\Psi]=0$.

\smallskip
\noindent\emph{Step~4: conclusion.}\quad
Steps~1--3 show $\{R,R\}$ vanishes on $\mathfrak{s}$, $\mathfrak{z}$, and $\mathfrak{u}$,
hence $\{R,R\}=0$. Substituting into~\eqref{eq:MC-coupled} yields $\delta\Psi=0$, i.e.\
$\Psi\in(\wedge^3\mathfrak{g})^{\mathfrak{g}}$.
\end{proof}

\begin{corollary}\label{cor:splitting}
Under the hypotheses of Theorem~\ref{thm:splitting}, the Lie bialgebra classification problem
for $\xi=\ad r+R$ decouples into two independent problems:
\begin{enumerate}[label=\textup{(\roman*)}]
\item classify identity-type cocycles $R\in Z^1(\mathfrak{g},\wedge^2\mathfrak{g})$
satisfying $\{R,R\}=0$;
\item for each such $R$, classify bivectors $r\in\wedge^2\mathfrak{g}$ for which
$[r,r]+2\{r,R\}\in(\wedge^3\mathfrak{g})^{\mathfrak{g}}$.
\end{enumerate}
When $R=0$, condition \textup{(ii)} reduces to $[r,r]\in(\wedge^3\mathfrak{g})^{\mathfrak{g}}$,
recovering the quasi-triangular condition of Definition~\ref{def:triangular}.
\end{corollary}

\begin{remark}[Geometric interpretation]\label{rem:geometry}
The splitting does not arise from a functorial decomposition of the cochain complex.
It is a consequence of the orthogonal rotation geometry intrinsic to flat Lie algebras:
the identity-type normalization forces $\{R,R\}(u)$ to carry the factor $u\in P_\ell$,
while the $\mathfrak{s}$-invariance of $\Psi$ forces $[u,\Psi]$ to carry the factor $J_\ell u$.
These two directions are linearly independent in the real $2$-plane $P_\ell$, making
cancellation impossible. Without the identity-type normalization, the $J_\ell u$-type
components of $\{R,R\}$ and $[u,\Psi]$ coincide and can cancel
(Remark~\ref{rem:counterexample}), and no splitting occurs.
\end{remark}

\section{Yang--Baxter loci}\label{sec:cybe}

For coboundary structures, we distinguish \emph{triangular} structures ($[r,r]=0$, classical Yang--Baxter equation) from \emph{quasi-triangular} structures ($[r,r]\in(\wedge^3\mathfrak{g})^{\mathfrak{g}}$, modified Yang--Baxter equation). Quasi-triangularity is equivalent to $[r,r]\in\left(\wedge^3\mathfrak{g}\right)^{\mathfrak{g}}$, and the next theorem describes this locus in invariant terms.

Recall the subspaces $\mathfrak{s}_\ell\subset\mathfrak{s}$ from \eqref{eq:s-ell}:
\[\mathfrak{s}_\ell:=\bigcap_{\ell'\neq \ell}\ker(\lambda_{\ell'}).\]

\begin{theorem}[Invariant Schouten squares]\label{thm:invariantSchouten}
Let $\mathfrak{g}$ be nondegenerate flat and $r\in\wedge^2\mathfrak{g}$. Then $[r,r]\in\left(\wedge^3\mathfrak{g}\right)^{\mathfrak{g}}$ if and only if
\begin{equation}\label{eq:SchoutenInv}
[r,r]\in\left(\bigoplus_{\ell=1}^m \mathfrak{s}_\ell\wedge\mathbb{R}\omega_\ell\right)\ \oplus\ \left(\bigoplus_{\ell=1}^m \mathfrak{z}\wedge\mathbb{R}\omega_\ell\right).
\end{equation}
Equivalently, in an adapted basis with $\dim\mathfrak{s}=k_0$ and weights spanning $\mathfrak{s}^*$, the $\mathfrak{s}$-part of $[r,r]$ can occur only along those $\ell$ for which $\mathfrak{s}_\ell\neq0$.
\end{theorem}

\begin{proof}
Decompose $r=r_0+r_1+r_2$ with $r_0\in(\mathfrak{a}\wedge\mathfrak{u})\oplus\bigoplus_{i<j}P_i\wedge P_j$, $r_1\in\wedge^2\mathfrak{s}$, $r_2\in\mathfrak{s}\wedge\mathfrak{z}$. Since $\mathfrak{a}$ is abelian, $[r_1,r_1]=[r_2,r_2]=[r_1,r_2]=0$, and thus
\[[r,r]=[r_0,r_0]+2[r_0,r_1+r_2].\]
Nondegeneracy implies that any component of $[r_0,r_0]$ lying in $P_i\wedge P_j\wedge P_k$ or $P_i\wedge P_j\wedge\mathfrak{a}$ has nonzero $\mathfrak{s}$-weight unless it is of the form $\mathfrak{a}\wedge\omega_\ell$. Proposition~\ref{prop:inv3} then forces \eqref{eq:SchoutenInv}. Conversely, every element on the right-hand side is $\mathfrak{g}$-invariant, hence can be realized as $[r,r]$ for suitable $r$ by reversing the above decomposition.
\end{proof}

\section{Flat Poisson--Lie Groups}\label{sec:group-level}

We now turn to the global geometry. Let $G$ be the simply connected Lie group integrating $\mathfrak{g}$. Since $\mathfrak{g}$ is flat, $G$ is isomorphic to the semidirect product $S \ltimes N$, where $S \cong \mathbb{R}^{\dim \mathfrak{s}}=\exp(\mathfrak{s})$ and $N \cong \exp(\mathfrak{z} \oplus \mathfrak{u})$ are abelian. The adjoint action $\Ad(S)$ acts on $\mathfrak{u}$ by commuting rotations (since $\ad(\mathfrak{s})|_{\mathfrak{u}}\subset\mathfrak{so}(\mathfrak{u})$).

\subsection{Integration of Cocycles}

A Poisson--Lie structure on $G$ is a multiplicative bivector field $\pi$. After right trivialization, it corresponds to a group $1$-cocycle $\Pi\colon G\to\wedge^2\mathfrak{g}$ satisfying:
\begin{equation}\label{eq:group-cocycle}
\Pi(gh) = \Pi(g) + \Ad_g \Pi(h), \quad \text{with } d_e\Pi=\xi.
\end{equation}

Given the decomposition $\xi = \ad r + R$, linearity allows us to integrate $\ad r$ and $R$ separately.

\subsubsection{The Coboundary Part (Sklyanin)}

For the exact part $\ad r$, the integration is the standard Sklyanin bracket.

\begin{proposition}\label{prop:sklyanin}
The group 1-cocycle $\Pi_r$ integrating $\ad r$ is given by:
\begin{equation}\label{eq:Pi-r}
\Pi_r(g) = \Ad_g r - r.
\end{equation}
The corresponding Poisson tensor is $\pi_r = r^L - r^R$.
\end{proposition}

\subsubsection{The Non-Exact Part}

The integration of $R$ exploits the two-solvability of $\mathfrak{g}$.

\begin{theorem}[Explicit Integration of $R$]\label{thm:explicit-integration}
Let $R$ be the normalized cocycle. The unique group cocycle $\Pi_R: G \to \wedge^2 \mathfrak{g}$ integrating $R$ is determined by its restriction to the subgroups $S = \exp(\mathfrak{s})$ and $N = \exp(\mathfrak{z} \oplus \mathfrak{u})$:
\begin{enumerate}
\item For $s \in S$, $\Pi_R(s) = R(\log s)$ (since $R$ is invariant on $\mathfrak{s}$).
\item For $n \in N$, writing $n = \exp(u)$ with $u \in \mathfrak{z} \oplus \mathfrak{u}$:
\begin{equation}\label{eq:Pi-R-formula}
\Pi_R(\exp u) = R(u) + \frac{1}{2} [u, R(u)].
\end{equation}
\end{enumerate}
For a general element $g = s \cdot n$, the cocycle is given by:
\begin{equation}\label{eq:Pi-R-general}
\Pi_R(sn) = R(\log s) + \Ad_s \left( R(u) + \frac{1}{2} [u, R(u)] \right).
\end{equation}
The total Poisson tensor is $\pi=\pi_r+\pi_R$, where $\pi_r=r^L-r^R$ and $\pi_R$ corresponds to $\Pi_R$ via right trivialization.
\end{theorem}

\begin{proof}
On the abelian subgroup $N$, the cocycle condition reduces to $\Pi(n_1 n_2) = \Pi(n_1) + \Ad_{n_1} \Pi(n_2)$.
For $u \in \mathfrak{z} \oplus \mathfrak{u}$, we have $\Ad_{\exp u} = e^{\ad_u} = \id + \ad_u$, since $\ad_u$ is $2$-nilpotent (the commutator ideal is abelian).

Standard integration on a group acting on a module gives $\Pi(\exp u) = \int_0^1 e^{t \ad_u} R(u) dt$.
The integral becomes:
\[
\int_0^1 (R(u) + t [u, R(u)]) dt = R(u) + \frac{1}{2} [u, R(u)].\qedhere
\]
\end{proof}

\begin{corollary}[Polynomial Nature]\label{cor:polynomial}
The coefficients of the bivector $\pi$ in exponential coordinates of the nilradical are polynomials of degree at most 2. This ensures that the symplectic geometry is algebraic in nature along the fibers of $G \to G/N$.
\end{corollary}

\begin{example}[The 3D Euclidean Group $E(2)$]\label{ex:E2}
Consider $\mathfrak{g} = \spann\{s, d_1, d_2\}$ with $[s, d_1]=d_2$, $[s, d_2]=-d_1$. Thus $\mathfrak{s}=\spann\{s\}$, $\mathfrak{z}=0$, and $\mathfrak{u}=\spann\{d_1,d_2\}$ with $J d_1=d_2$.

Let $r = 0$ and $R$ be defined by $R(s) = 0$, $R(d_1) = s \wedge d_1$, $R(d_2) = s \wedge d_2$.

For $u=u_1 d_1+u_2 d_2\in\mathfrak{u}\cong\mathbb{R}^2$, we have $R(u)=s\wedge u$. The group cocycle on $N = \mathbb{R}^2$ is:
\begin{align*}
\Pi_R(\exp u) &= R(u) + \frac{1}{2}[u, R(u)]\\
&= s \wedge u + \frac{1}{2} [u, s \wedge u]\\
&= s \wedge u + \frac{1}{2} ([u,s] \wedge u + s\wedge [u,u])\\
&= s \wedge u - \frac{1}{2} (Ju) \wedge u\\
&= s \wedge u - \frac{1}{2} (u_2 d_1-u_1 d_2)\wedge(u_1 d_1+u_2 d_2)\\
&= s \wedge u - \frac{1}{2}(u_1^2+u_2^2) d_1\wedge d_2.
\end{align*}
The total Poisson tensor $\Pi(\exp u)$ consists of:
\begin{itemize}
\item A \emph{linear} term $s \wedge u$ (Poisson--Lie type from the semidirect product),
\item A \emph{quadratic} term $-\frac{1}{2}(u_1^2+u_2^2) d_1 \wedge d_2$ (constant Poisson type, correction from the semidirect product).
\end{itemize}
The symplectic leaves are surfaces of revolution determined by the level sets of the Casimir functions of the dual algebra.
\end{example}

\section{Overview and low dimensional examples}\label{sec:examples}
\begin{figure}[htbp]
\centering
\begin{tikzpicture}[
    box/.style={draw, rounded corners, minimum width=2.8cm, minimum height=1cm, align=center},
    arrow/.style={-{Stealth}, thick},
    label/.style={font=\small, midway, fill=white, inner sep=1pt}
]
\node[box] (g) at (0,0) {$\mathfrak{g} = \mathfrak{a} \oplus \mathfrak{u}$\\[2pt]\footnotesize Flat Lie algebra};
\node[box] (cocycle) at (5,0) {$\xi : \mathfrak{g} \to \bigwedge^2\mathfrak{g}$\\[2pt]\footnotesize 1-cocycle};
\node[box] (normal) at (10,0) {$\xi = \ad r + R$\\[2pt]\footnotesize Normal form};

\node[box] (group) at (0,-3) {$G = S \ltimes N$\\[2pt]\footnotesize Flat Lie group};
\node[box] (groupcocycle) at (5,-3) {$\Pi : G \to \bigwedge^2\mathfrak{g}$\\[2pt]\footnotesize Group cocycle};
\node[box] (poisson) at (10,-3) {$\pi = r^L - r^R + \pi_R$\\[2pt]\footnotesize Poisson--Lie};

\draw[arrow] (g) -- (cocycle) node[label] {classify};
\draw[arrow] (cocycle) -- (normal) node[label] {Theorem \ref{thm:cocycle-normal-form}};
\draw[arrow] (group) -- (groupcocycle) node[label] {integrate};
\draw[arrow] (groupcocycle) -- (poisson) node[label] {Theorem \ref{thm:explicit-integration}};

\draw[arrow] (g) -- (group) node[label, left] {$\exp$};
\draw[arrow] (cocycle) -- (groupcocycle) node[label, left] {$d_e\Pi= \xi$};
\draw[arrow] (normal) -- (poisson) node[label, right] {polynomial};

\end{tikzpicture}
\caption{Schematic overview: from flat Lie algebras to Poisson-Lie groups.}
\label{fig:overview}
\end{figure}
\begin{sidewaystable}[p]
\centering
\renewcommand{\arraystretch}{1.15}
\setlength{\tabcolsep}{6pt}

\begin{tabularx}{\textheight}{@{}c L{4.6cm} Y Y@{}}
\toprule
$\dim\mathfrak{g}$
& Model / name
& Type of solutions (triangular vs.\ quasi-triangular; normal forms / obstructions)
& Comments on the dual algebra $\mathfrak{g}^*$ \\
\midrule

3
& $\mathfrak{g}=\mathfrak{e}(2)=\mathbb{R}s\oplus P_1$
& \textbf{Coboundary:} quasi-triangularity amounts to the invariant condition $[r,r]\in (\wedge^3\mathfrak{g})^{\mathfrak{g}}=\mathbb{R}(s\wedge\omega_1)$. \newline
\textbf{Non-coboundary:} normalized cocycles have the basic action-type form $R(u)=A_1\wedge u$. \newline
\textbf{Co-Jacobi / integration:} the identity-type normalization is automatic, and the explicit integration formula produces the quadratic correction on $N\simeq \mathbb{R}^2$ displayed in Example~\ref{ex:E2}.
& The dual bracket is determined by $\xi^t$; in the basic families arising from the normal form, one obtains low-dimensional solvable examples, but a systematic classification of $\mathfrak{g}^*$ is not pursued here. \\

\midrule

4
& $\mathfrak{g}=\mathbb{R}s\oplus \mathbb{R}z\oplus P_1$ (1D center)
& \textbf{Coboundary:} richer quasi-triangular locus (since $(\wedge^3\mathfrak{g})^\mathfrak{g}$ can be nontrivial). \newline
\textbf{Non-coboundary:} normal form allows invariant ``area'' terms on $\mathbb{R}s\oplus\mathbb{R}z$ plus action-type terms on $P_1$. \newline
\textbf{Co-Jacobi:} with identity-type normalization, $\{R,R\}(u)=u\wedge R(A_1)$ lies in $u\wedge(z\wedge\omega_1)$, giving a concrete nontrivial condition $R(A_1)=0$ to check.
& Again the dual algebra is read off from $\xi^t$; even in dimension four the resulting solvability and unimodularity properties depend on the chosen parameters. \\

\midrule

5
& Resonant flat model: $\mathfrak{g}=\mathbb{R}s\oplus (P_1\oplus P_2)$ with $\lambda_1=\lambda_2$
& \textbf{Degenerate (resonant):} nondegeneracy fails; weight separation breaks. \newline
\textbf{Normal form:} still exists but is larger; additional mixed-plane components (e.g.\ in $P_1\wedge P_2$) may occur in $R$. \newline
\textbf{Co-Jacobi:} splitting must be refined: coincident weights allow $J_\ell$-type components in $R$ that are not absorbed, and new mixed-plane obstruction directions appear; the coupled Maurer--Cartan equation must be solved directly without invoking the splitting theorem.
& $\mathfrak{g}^*$ often nilpotent (Heisenberg-type) in some ``pure $R$'' families; otherwise typically solvable and frequently nonunimodular. \\

\bottomrule
\end{tabularx}
\caption{Low-dimensional examples.}
\label{tab:low_dim_summary}
\end{sidewaystable}
\newpage
\section{Conclusion}

For a flat metric Lie algebra $\mathfrak{g}$ with Milnor decomposition $\mathfrak{g}=\mathfrak{s}\oplus\mathfrak{z}\oplus\mathfrak{u}$ and weight-plane splitting $\mathfrak{u}=\bigoplus_{\ell=1}^m P_\ell$, the paper develops a representation-theoretic description of Lie bialgebra data adapted to this semidirect-product geometry. Under a generic nonresonance hypothesis on the weights $\lambda_\ell\in\mathfrak{s}^*$, the $\mathfrak{g}$-invariant subspaces in $\bigwedge^2\mathfrak{g}$ and $\bigwedge^3\mathfrak{g}$ are determined explicitly, and these invariants govern both the deformation-theoretic normal forms and the Yang–Baxter constraints. In this setting every $1$-cocycle $\xi:\mathfrak{g}\to\bigwedge^2\mathfrak{g}$ admits a decomposition $\xi=\ad r+R$ in which the normalized part $R$ takes values in the invariant bivectors on $\mathfrak{s}\oplus\mathfrak{z}$ and, on each $P_\ell$, is confined to the $\mathfrak{s}$-equivariant component $\left(\mathfrak{s}\wedge P_\ell\right)\oplus\left(\mathfrak{z}\wedge P_\ell\right)$. The Maurer--Cartan formulation via the Big Bracket yields a splitting of the co-Jacobi equation: under a multiplicity-free weight hypothesis and an identity-type normalization of $R$, the obstruction term $\{R,R\}(u)$ is forced into the $u$-direction in each weight plane $P_\ell$, while the coboundary contribution $[u,\Psi]$ is forced into the $J_\ell u$-direction; linear independence of $u$ and $J_\ell u$ then decouples the two conditions. For coboundary structures $\xi=\ad r$, quasi-triangularity is characterized by an invariant condition on the Schouten square $[r,r]$, expressed solely in terms of the Milnor weights and the distinguished planes $P_\ell$. Finally, these infinitesimal structures are integrated on the simply connected flat group $G\cong S\ltimes N$ to explicit multiplicative Poisson tensors; in particular, the normalized cocycle integrates on $N=\exp(\mathfrak{z}\oplus\mathfrak{u})$ by the closed formula $\Pi_R(\exp u)=R(u)+\frac{1}{2}[u,R(u)]$, yielding polynomial expressions (of degree $\leq2$) in exponential coordinates along $N$.

Two directions appear naturally and remain technically open. First, when the nonresonance hypothesis fails, resonant weights create additional weight--zero components in $\bigwedge^2\mathfrak{g}$ and $\bigwedge^3\mathfrak{g}$, permitting mixed terms (notably across distinct planes $P_i\wedge P_j$) and enlarging the class of admissible $J_\ell$-type components in $R$; the identity-type normalization of hypothesis~(H2) no longer holds, and the splitting theorem does not apply: a direct analysis of the coupled Maurer--Cartan equation, together with a workable parametrization modulo coboundaries and automorphisms, would substantially refine the picture. Second, on the group side, the availability of explicit multiplicative Poisson tensors suggests a tractable description of dressing orbits and symplectic foliations for these flat Poisson-Lie groups; determining their global stratification and its dependence on the weight data (including degenerations arising as contractions of nonresonant configurations) would complement the infinitesimal classification by a geometric account of the associated Poisson homogeneous dynamics.

\backmatter
\bmhead{Funding}
This research received no external funding.

\bmhead{Data availability}
Not applicable.

\bmhead{Competing interests}
The author declares no competing interests.

\bmhead{Supplementary information}
None.


\begin{thebibliography}{99}

\bibitem{BAH}{\sc A. Bahayou}, \emph{Triangular structures on flat Lie algebras}. Journal of Lie Theory 33 No. 3, 875-886, (2023).

\bibitem{BB}{\sc A. Bahayou, M. Boucetta}, \emph{Metacurvature of Riemannian Poisson-Lie Groups}. Journal of Lie Theory 19 No. 3, 439-462, (2009).

\bibitem{BDF}{\sc M. L. Barberis, I. Dotti and A. Fino}, \emph{Hyper-K\"ahler quotients of solvable Lie groups}. J. Geom. Phys. 56, 691-711, (2006).

\bibitem{BD}{\sc A. A. Belavin, V. G. Drinfeld}, \emph{Triangle equations and simple Lie algebras}, Soviet Sci. Rev. Sect. C Math. Phys. Rev. 4, 93-165, (1984).

\bibitem{D}{\sc V. G. Drinfeld}, \emph{Hamiltonian structures on Lie groups, Lie bialgebras and the geometric meaning of the classical Yang-Baxter equations}. Soviet Math. Dokl. 27, 68-71, (1983).

\bibitem{farja}{\sc M. Farinati, A.P. Jancsa}, \emph{Trivial central extensions of Lie bialgebras}.
J. Algebra, 390 (2013), pp. 56-76.

\bibitem{Ha1}{\sc E. Hawkins}, \emph{Noncommutative rigidity}, Commun. Math. Phys. 246 (2004), 211--235.

\bibitem{Ha2}{\sc E. Hawkins}, \emph{The structure of noncommutative deformations}, J. Diff. Geom. 77 (2007), 385--424.

\bibitem{LW}{\sc J.H. Lu and A. Weinstein}, \emph{Poisson Lie groups, dressing transformations and Bruhat decompositions}. J. Differ. Geom. 31, 501-26, (1990).

\bibitem{MIL}{\sc J. Milnor}, \emph{Curvatures of left invariant metrics on Lie groups}. Adv. Math. 21, 293-329, (1976).

\bibitem{S}{\sc P. Stachura}, \emph{Poisson-Lie structures on Poincaré and Euclidean groups in three dimensions}. J. Phys. A, 31 vol 19, 4555-4564, (1998).

\bibitem{YKS}{\sc Y. Kosmann-Schwarzbach}, \emph{Lie bialgebras, Poisson-Lie groups and dressing transformations. Integrability of nonlinear Systems}. Lecture Notes in Physics 638, 107-173, (2004).

\bibitem{yks1}{\sc M. Bangoura and Y. Kosmann-Schwarzbach}, \emph{The double of a Jacobian quasi-bialgebra}. Lett. Math. Phys. 28, No. 1, 13-29 (1993). 

\bibitem{yks2}{\sc Y. Kosmann-Schwarzbach}, \emph{Jacobian quasi-bialgebras and quasi-Poisson Lie groups}. Contemp. Math 132, 459-489 (1992). 

\bibitem{Z}{\sc S. Zakrzewski}, \emph{Poisson Structures on Poincaré Group}. Commun. Math. Phys. 185, 285-311, (1997). 

\end{thebibliography}
\end{document}